# Analytic solution of linear fractional differential equation with Jumarie derivative in term of Mittag-Leffler function

Uttam Ghosh (1), Srijan Sengupta (2a), Susmita Sarkar (2b), Shantanu Das (3)
(1): Department of Mathematics, Nabadwip Vidyasagar College, Nabadwip, Nadia, West Bengal, India; Email: uttam_math@yahoo.co.in
(2):Department of Applied Mathematics, University of Calcutta, Kolkata, India
Email (2b): susmita62@yahoo.co.in,
(3)Scientist H+, RCSDS, Reactor Control Div. (complex) BARC Mumbai India
Senior Research Professor, Dept. of Physics, Jadavpur University Kolkata
Adjunct Professor. DIAT-Pune
Ex- UGC Visiting Fellow. Dept of Appl. Mathematics; Univ. of Calcutta
Email (3): shantanu@barc.gov.in

**Abstract**

There is no unified method to solve the fractional differential equation. The type of derivative here used in this paper is of Jumarie formulation, for the several differential equations studied. Here we develop an algorithm to solve the linear fractional differential equation composed via Jumarie fractional derivative in terms of Mittag-Leffler function; and show its conjugation with ordinary calculus. In these fractional differential equations the one parameter Mittag-Leffler function plays the role similar as exponential function used in ordinary differential equations.

**1.0 Introduction**

The analytical solutions of the fractional differential equation are emerging branch of applied science also in basic science. Different methods are developing to solve the fractional differential equations. Since the definition of fractional derivative is modifying to relate it with the classical derivative. Mathematicians are trying to develop the formulas of fractional calculus but geometry of fractional derivative has no concrete shape [1]. Depending on different type of derivatives different methods of solution are developing [2-7]. Riemann-Liouville definition the fractional derivative of a constant is non-zero which creates a difficulty to relate between the classical calculus. To overcome this difficulty Jumarie [2]-[5] modified the definition of fractional derivative of Riemann-Liouvell type and with this new formulation , we obtain the derivative of a constant as zero. Thus using this definition linking between the fractional and classical calculus becomes easier. There is no unique method to solve the linear fractional differential equations. Using the Jumarie modified definition of fractional derivative we obtain the derivative of Mittag-Leffler function as Mittal-Leffler function, as in case of classical whole-order derivative the derivative of $\exp(x)$ is itself exponential function. Thus via use of Jumarie modified Riemann-Liouvelli derivative, there exists conjugation with classical calculus, which eases in many cases to solve fractional differential equation composed with Jumarie fractional derivative. Here we want to develop an algorithm to solving the linear fractional differential equation using the Mittag-Leffler function. We have obtained applied this method to homogeneous fractional differential equations and got corresponding fundamental solution.

Organization of the paper is as follows. In section 2.0 some definition of fractional derivative is reproduced with essential examples. In section 3.0 and 4.0 some properties of Mittag-Leffler function is described. Finally in section 5.0 the methods for solving the linear fractional differential equation composed by Jumarie fractional derivative is developed using the Mittag-Leffler function.

**2.0 Some definitions of fractional**

There are many definition of fractional derivative. Grunwald-Letnikov fractional derivative [6], Liouville fractional derivative [8], Riemann-Liouville fractional derivative [10], Caputo fractional derivative [8-11], Kolwanker-Gangal local fractional derivative [12-16], Jumarie modified fractional derivative [2]. Here we use the Riemann-Liouville fractional derivative and its modified form by Jumarie [2].

**2.1 Riemann-Liouville definition of fractional derivative**

Let the function $f(t)$ is one time integrable then the integro-differential expression

$$_aD_t^\alpha f(t) = \frac{1}{\Gamma(-\alpha+m+1)} \left(\frac{d}{dt}\right)^{m+1} \int_a^t (t-\tau)^{m-\alpha} f(\tau) d\tau$$

is known as the Riemann-Liouville (R-L) definition of fractional derivative [6] with $m$ as integer with condition $m \leq \alpha < m+1$.

In Riemann-Liouville definition the function, $f(t)$ is integrated $(m-\alpha)$ fold and then differentiate $m+1$ times. We can re-write the above as follows

The left R-L fractional derivative is defined by

$$_aD_t^\alpha f(t) = \frac{1}{\Gamma(k-\alpha)} \left(\frac{d}{dt}\right)^k \int_a^t (t-\tau)^{k-\alpha-1} f(\tau) d\tau$$

And the right R-L derivative is

$$_tD_b^\alpha f(t) = \frac{1}{\Gamma(k-\alpha)} \left(-\frac{d}{dt}\right)^k \int_t^b (\tau-t)^{k-\alpha-1} f(\tau) d\tau$$

Where in above $k$ is integer such that $(k-1) \leq \alpha < k$ that is just greater than fractional number $\alpha$.

Using the left R-L derivative we get the fractional derivative of the function $f(t) = K$ as non-zero, as demonstrated below.

$$_aD_t^\alpha f(t) = \frac{1}{\Gamma(1-\alpha)} \frac{d}{dt} \int_a^t (t-\xi)^{-\alpha} K d\xi$$

$$= -\frac{K}{\Gamma(1-\alpha)} \frac{d}{dt} \frac{(t-\xi)^{1-\alpha}}{1-\alpha}\bigg]_a^t$$

$$= K \frac{(t-a)^{1-\alpha}}{\Gamma(1-\alpha)}$$

Similarly the right R-L derivative of $f(t) = K$ is

$$_tD_b^\alpha f(t) = K \frac{(b-t)^{1-\alpha}}{\Gamma(1-\alpha)}.$$

This shows that the fractional derivative of a constant ($K$) is non-zero but in classical calculus derivative of a constant is zero which is contradiction between the classical derivative and the fractional derivative of a constant. To overcome this difference Jumarie [2] modified the left R-L fractional derivative.

We get the R-L left derivative of a power function as

**2. 2 Jumarie modified definition of the fractional derivative is**

$$D_t^\alpha f(t) \stackrel{\text{def}}{=} \begin{cases} \frac{1}{\Gamma(-\alpha)} \int_0^t (t-\xi)^{-\alpha-1} f(\xi) d\xi, & \alpha < 0. \\ \frac{1}{\Gamma(1-\alpha)} \frac{d}{dt} \int_0^t (t-\xi)^{-\alpha} [f(\xi) - f(0)] d\xi, & 0 < \alpha < 1. \\ [f^{(\alpha-n)}(t)]^{(n)}, & n \leq \alpha < n+1, \quad n \geq 1. \end{cases}$$

Using this definition we get $D^\alpha\{K\} = 0$, $\quad 0 \leq \alpha < 1$.

The above formula in line-1, becomes fractional order integration if we replace $\alpha$ by $-\alpha$ which is

$$_aD_t^{-\alpha} f(t) = \frac{1}{\Gamma(\alpha)} \int_a^t (t-\tau)^{\alpha-1} f(\tau) d\tau \tag{2}$$

Using the above formula we get for $f(t) = (t-a)^\gamma$, the fractional integral for order $\alpha$

$$_aD_t^{-\alpha}(t-a)^\gamma = \frac{1}{\Gamma(\alpha)} \int_a^t (t-\tau)^{\alpha-1} (\tau-a)^\gamma d\tau$$

Using the substitution $\tau = a + \xi(t-a)$ we have for; $\tau = a$, $\xi = 0$ and for $\tau = t, \xi = 1$; $d\tau = (t-a)d\xi$, $(t-\tau) = t - a - \xi(t-a) = (t-a)(1-\xi)$; $(\tau - a) = \xi(t-a)$, we get the following

$$_aD_t^{-\alpha}(t-a)^\gamma = \frac{1}{\Gamma(\alpha)} \int_a^t (t-\tau)^{\alpha-1}(\tau-a)^\gamma d\tau$$

$$= \frac{1}{\Gamma(\alpha)} \int_0^1 (t-a)^{\alpha-1}(1-\xi)^{\alpha-1}\xi^\gamma (t-a)^\gamma (t-a)d\xi$$

$$= \frac{(t-a)^{\gamma+\alpha}}{\Gamma(\alpha)} \int_0^1 \xi^\gamma (1-\xi)^{\alpha-1} d\xi$$

$$= \frac{(t-a)^{\gamma+\alpha}}{\Gamma(\alpha)} B(\alpha, \gamma+1)$$

$$= \frac{\Gamma(\gamma+1)}{\Gamma(\gamma+1+\alpha)} (t-a)^{\gamma+\alpha}, \qquad (\alpha < 0, \gamma > -1)$$

We used Beta-function $B(\alpha, \gamma+1) = \int_0^1 \xi^\gamma (1-\xi)^{\alpha-1} d\xi = \frac{\Gamma(\alpha)\Gamma(\gamma+1)}{\Gamma(\alpha+\gamma+1)}$ defined as

$$B(p,q) \stackrel{\text{def}}{=} \int_0^1 u^{p-1}(1-u)^{q-1} du = \frac{\Gamma(p)\Gamma(q)}{\Gamma(p+q)}$$

Applying the above obtained result the fractional integral of order $(1-\upsilon)$, with $0 \leq \upsilon < 1$ is

$$_aD_t^{-(1-\upsilon)}(t-a)^\gamma = \frac{\Gamma(\gamma+1)}{\Gamma(\gamma+2-\upsilon)} (t-a)^{\gamma+1-\upsilon}$$

Taking one whole derivative of the above we get

$$D^1 \left[ _aD_t^{-(1-\upsilon)}(t-a)^\gamma \right] = {}_aD_t^\upsilon (t-a)^\gamma$$

$$= \frac{d}{dt} \left[ \frac{\Gamma(\gamma+1)}{\Gamma(\gamma+2-\upsilon)} (t-a)^{\gamma+1-\upsilon} \right]$$

$$= \frac{\Gamma(\gamma+1)}{\Gamma(\gamma-\upsilon+2)} (\gamma+1-\upsilon)(t-a)^{\gamma-\upsilon}$$

$$= \frac{\Gamma(\gamma+1)}{(\gamma-\upsilon+1)\Gamma(\gamma+1-\upsilon)} (\gamma+1-\upsilon)(t-a)^{\gamma-\upsilon}$$

$$= \frac{\Gamma(\gamma+1)}{\Gamma(\gamma+1-\upsilon)} (t-a)^{\gamma-\upsilon}$$

We have therefore calculated fractional derivative by R-L left formula for $\upsilon$ such that $0 \leq \upsilon < 1$ thus our nearest integer is one that is $k = 1$ and we write that below

$$_aD_t^v f(t) = \frac{1}{\Gamma(1-v)}\left(\frac{d}{dt}\right)\int_a^t (t-\tau)^{-v} f(\tau)d\tau$$

$$= D^1\left[D_t^{-(1-v)} f(t)\right]$$

Thus for $a=0$, the fractional RL derivative of $f(t)=t^\gamma$ is

$$_0D_t^v t^\gamma = \frac{\Gamma(\gamma+1)}{\Gamma(\gamma+1-v)} t^{\gamma-v}$$

For a constant function $f(t)=1$, putting in above expression $\gamma=0$, we get

$$_0D_t^v[1] = \frac{1}{\Gamma(1-v)} t^{-v}$$

Let us now see what Jumarrie derivative is from above R-L derivative obtained for $f(t)=t^\gamma$. The composition of the Jumarie derivative, with start point of integration as $t=a$ and $f(a)=a^\gamma$ is

$$f^{(\alpha)}[t^\gamma]_a^t = \frac{1}{\Gamma(1-\alpha)} \frac{d}{dt} \int_a^t (t-\xi)^{-\alpha}\left[\xi^\gamma - a^\gamma\right] d\xi$$

$$D^1\left[_aD_t^{-(1-\alpha)} t^\gamma - {_aD_t^{-(1-\alpha)}} a^\gamma\right] = {_aD_t^\alpha} t^\gamma - {_aD_t^\alpha} a^\gamma$$

The above expression show that this Jumarrie derivative is composed of two RL derivatives those are $_aD_t^\alpha t^\gamma$ minus RL derivative of a constant $_aD_t^\alpha a^\gamma$. Going by similar steps as done for $_0D_t^\alpha t^\gamma$, we get first the fractional integral in terms of incomplete Gamma function as

$$_aD_t^{-(1-\alpha)} t^\gamma = \frac{t^{\gamma+1-\alpha}}{\Gamma(1-\alpha)} \int_0^{(t-a)/t} z^{-\alpha}(1-z)^\gamma dz = \frac{t^{\gamma+1-\alpha}}{\Gamma(1-\alpha)} B_\eta(1-\alpha, \gamma+1) \qquad \eta = \frac{t-a}{t}$$

The fractional derivative of $t^\gamma$ is by taking one whole derivative of above expression we get the following

$$_aD_t^\alpha t^\gamma = \frac{d}{dt}\left[\frac{t^{\gamma+1-\alpha}}{\Gamma(1-\alpha)} B_\eta(1-\alpha, \gamma+1)\right] \qquad \eta = \frac{t-a}{t}$$

The fractional derivative of $a^\gamma$ is

$$_aD_t^\alpha a^\gamma = \frac{a^\gamma}{\Gamma(1-\alpha)} (t-a)^{-\alpha}$$

Therefore

$$f^{(\alpha)}[t^\gamma]_a^t = \frac{1}{\Gamma(1-\alpha)} \frac{d}{dt} \int_a^t (t-\xi)^{-\alpha} \left[\xi^\gamma - a^\gamma\right] d\xi$$

$$= D^1\left[{}_aD_t^{-(1-\alpha)} t^\gamma - {}_aD_t^{-(1-\alpha)} a^\gamma\right] = {}_aD_t^\alpha t^\gamma - {}_aD_t^\alpha a^\gamma$$

$$= \frac{d}{dt}\left[\frac{t^{\gamma+1-\alpha}}{\Gamma(1-\alpha)} B_\eta(1-\alpha, \gamma+1)\right] - \frac{a^\gamma}{\Gamma(1-\alpha)}(t-a)^{-\alpha}$$

For $a=0$, we have

$$f^{(\alpha)}[t^\gamma]_0^t = \frac{1}{\Gamma(1-\alpha)} \frac{d}{dt} \int_0^t (t-\xi)^{-\alpha}\left[\xi^\gamma - 0\right] d\xi$$

$$= D^1\left[{}_0D_t^{-(1-\alpha)} t^\gamma\right] = {}_0D_t^\alpha t^\gamma = \frac{\Gamma(\gamma+1)}{\Gamma(\gamma+1-\alpha)}(t)^{\gamma-\alpha}$$

We will be using Jumarie derivative for power function $t^\gamma$ with start point of differentiation as $a=0$, in subsequent sections. When start point of differentiation is non-zero we will be shifting the origin to that non-zero point and use the above formula.

### 3.0 Some properties Mittag-Leffler function and its application

In 1903 Mittag-Leffler [17]-[19] introduce a function defined by an infinite series

$$E_\alpha(at^\alpha) \overset{\text{def}}{=} 1 + \frac{at^\alpha}{\Gamma(1+\alpha)} + \frac{a^2 t^{2\alpha}}{\Gamma(1+2\alpha)} + \frac{a^3 t^{3\alpha}}{\Gamma(1+3\alpha)} + \ldots$$

is the one parameter Mittag-Leffler function.

Using Jumarie derivative of order $\alpha$, with $0 \le \alpha < 1$ with start point as $a=0$ for $f(t) = t^{n\alpha}$, that is $(t)^{\alpha(n-1)}\Gamma(n\alpha+1)/\Gamma[\alpha(n-1)+1]$, for $n=1,2,3,\ldots$; and also using Jumarie derivative of constant as zero, we get the following very useful identity. In all the subsequent sections we will say $D^\alpha$ is the Jumarie derivative with zero as start point

$$D^\alpha(E_\alpha(at^\alpha)) = D^\alpha\left(1 + \frac{at^\alpha}{\Gamma(1+\alpha)} + \frac{a^2 t^{2\alpha}}{\Gamma(1+2\alpha)} + \frac{a^3 t^{3\alpha}}{\Gamma(1+3\alpha)} + \ldots\right)$$

$$= 0 + \frac{\Gamma(1+\alpha)a}{\Gamma(1)\Gamma(1+\alpha)} + \frac{\Gamma(1+2\alpha)a^2 t^\alpha}{\Gamma(1+2\alpha)\Gamma(1+\alpha)} + \frac{\Gamma(1+3\alpha)a^3 t^{2\alpha}}{\Gamma(1+3\alpha)\Gamma(1+2\alpha)} + \ldots$$

$$= a\left(1 + \frac{at^\alpha}{\Gamma(1+\alpha)} + \frac{a^2 t^{2\alpha}}{\Gamma(1+2\alpha)} + \frac{a^3 t^{3\alpha}}{\Gamma(1+3\alpha)} + \ldots\right)$$

$$= aE_\alpha(at^\alpha)$$

Thus

$$D^{\alpha}(E_{\alpha}(at^{\alpha})) = aE_{\alpha}(at^{\alpha}) \qquad (1)$$

This shows that $AE_{\alpha}(at^{\alpha})$ is a solution is a solution of the fractional differential equation

$$D^{\alpha}y = ay \qquad (2)$$

Where $A$ is arbitrary constant.

Therefore

$$D^{\alpha}y = ay$$

with $y(0) = 1$ has solution

$$y = E_{\alpha}(at^{\alpha}).$$

Using this property of the Mittag-Leffler one can easily prove the following theorem.

**Theorem 1:** The Mittag-Leffler function $E_{\alpha}(at^{\alpha})$ satisfies the relation

$$E_{\alpha}(at^{\alpha})E_{\alpha}(bt^{\alpha}) = E_{\alpha}((a+b)t^{\alpha})$$

Proof: Let $y = E_{\alpha}(at^{\alpha})E_{\alpha}(bt^{\alpha})$ then

$$\begin{aligned} y &= E_{\alpha}(at^{\alpha})E_{\alpha}(bt^{\alpha}) \\ D^{\alpha}y &= E_{\alpha}(bt^{\alpha})D^{\alpha}(at^{\alpha}) + E_{\alpha}(at^{\alpha})D^{\alpha}(bt^{\alpha}) \\ D^{\alpha}y &= aE_{\alpha}(at^{\alpha})E_{\alpha}(bt^{\alpha}) + bE_{\alpha}(at^{\alpha})E_{\alpha}(bt^{\alpha}) \\ &= (a+b)E_{\alpha}(at^{\alpha})E_{\alpha}(bt^{\alpha}) \\ D^{\alpha}y &= (a+b)y \end{aligned} \qquad (3)$$

Using the solution of the equation (2) we get the solution of the equation (3) in the following form

$$y = AE_{\alpha}\left([a+b]t^{\alpha}\right)$$

From the definition of $y$ we get $y(0) = 1$. Therefore we have $y = E_{\alpha}((a+b)t^{\alpha})$.

Thus we get

$$E_{\alpha}((a+b)t^{\alpha}) = E_{\alpha}(at^{\alpha})E_{\alpha}(bt^{\alpha}) \qquad (4)$$

We get useful property of one parameter Mittag-Leffler function.

Using the above property of Mittag-Leffler function we get

$$E_\alpha(at^\alpha)E_\alpha(-at^\alpha) = 1 \quad \text{or} \quad E_\alpha(-at^\alpha) = \frac{1}{E_\alpha(at^\alpha)}. \tag{5}$$

$$E_\alpha(at^\alpha)E_\alpha(at^\alpha) = E_\alpha(2at^\alpha) \tag{6}$$

**4.0 Complex Mittag-Leffler function and its properties**

Jumarie [2008] defined the complex Mittag-Leffler in the following form

$$E_\alpha(it^\alpha) \stackrel{def}{=} \cos_\alpha(t^\alpha) + i\sin_\alpha(t^\alpha)$$

$$\cos_\alpha(t^\alpha) = \frac{E_\alpha(it^\alpha) + E_\alpha(-it^\alpha)}{2} = \sum_{k=1}^{\infty}(-1)^k \frac{t^{2k\alpha}}{(2k\alpha)!}$$

$$\sin_\alpha(t^\alpha) = \frac{E_\alpha(it^\alpha) - E_\alpha(-it^\alpha)}{2} = \sum_{k=1}^{\infty}(-1)^k \frac{t^{(2k+1)\alpha}}{(2k\alpha+\alpha)!}$$

On the other hand Jumarie [2008] defined period ($M_\alpha$) of the function $E_\alpha(it^\alpha)$ in the following form, taking $E_\alpha(i(M_\alpha)^\alpha) = 1$ and therefore

$$\cos_\alpha(t+M_\alpha)^\alpha = \cos_\alpha(t^\alpha) \quad\quad \sin_\alpha(t+M_\alpha)^\alpha = \sin_\alpha(t^\alpha)$$
$$\cos_\alpha((-t)^\alpha) = \cos_\alpha(t^\alpha) \quad\quad \sin_\alpha((-t)^\alpha) = (-1)^\alpha \sin_\alpha(t^\alpha).$$

The series presentation of $\cos_\alpha(t^\alpha)$ is

$$\cos_\alpha(t^\alpha) = 1 - \frac{t^{2\alpha}}{\Gamma(1+2\alpha)} + \frac{t^{4\alpha}}{\Gamma(1+4\alpha)} - \frac{t^{6\alpha}}{\Gamma(1+6\alpha)} + \ldots$$

Taking term by term Jumarie derivative we get

$$D^\alpha[\cos_\alpha(t^\alpha)] = 0 - \frac{\Gamma(1+2\alpha)t^{2\alpha-\alpha}}{\Gamma(1+2\alpha)\Gamma(1+\alpha)} + \frac{\Gamma(1+4\alpha)t^{4\alpha-\alpha}}{\Gamma(1+4\alpha)\Gamma(1+3\alpha)} - \frac{\Gamma(1+6\alpha)t^{6\alpha-\alpha}}{\Gamma(1+6\alpha)\Gamma(1+5\alpha)} + \ldots$$

$$= -\left[\frac{t^\alpha}{\Gamma(1+\alpha)} - \frac{t^{3\alpha}}{\Gamma(1+3\alpha)} + \ldots\right] = -\sin_\alpha(t^\alpha)$$

**5.0 Solution of linear second order fractional differential equation**

Let us consider the function

$$y = AE_\alpha(at^\alpha) + BE_\alpha(bt^\alpha)$$

with *A* and *B* is constants. Differentiating $\alpha$–times with respect to *t*, for $0 < \alpha < 1$, with Jumarie derivative we get

$$D^\alpha y = AaE_\alpha(at^\alpha) + BbE_\alpha(bt^\alpha)$$
$$D^\alpha y - ay = AaE_\alpha(at^\alpha) + BbE_\alpha(bt^\alpha) - ay$$
$$= AaE_\alpha(at^\alpha) + BbE_\alpha(bt^\alpha) - a\left(AE_\alpha(at^\alpha) + BE_\alpha(bt^\alpha)\right)$$
$$= B(b-a)E_\alpha(bt^\alpha)$$
$$D^\alpha y - ay = B(b-a)E_\alpha(bt^\alpha)$$

Differentiating above by Jumarrie derivative and re-arranging, we get

$$D^{2\alpha} y - aD^\alpha y = Bb(b-a)E_\alpha(bt^\alpha)$$
$$D^{2\alpha} y - (a+b)D^\alpha y + aby = 0$$

This shows that the fractional differential equation

$$D^{2\alpha} y - (a+b)D^\alpha y + aby = 0$$

has solution in the form

$$y = AE_\alpha(at^\alpha) + BE_\alpha(bt^\alpha).$$

On the other hand consider the differential equation

$$D^{2\alpha} y - (a+b)D^\alpha y + aby = 0$$

it can be express in the following form

$$(D^\alpha - a)(D^\alpha - b)y(t) = 0. \tag{7}$$

Let, $(D^\alpha - b)y(t) = x(t)$ then equation (7) reduce to the form

$$(D^\alpha - a)x(t) = 0 \quad \text{or} \quad D^\alpha x(t) = ax(t)$$

Solution of the above equation is same as the solution of the equation (2) which is

$$x(t) = A_1 E_\alpha(at^\alpha)$$
$$(D^\alpha - b)y(t) = A_1 E_\alpha(at^\alpha)$$
$$D^\alpha y - by = A_1 E_\alpha(at^\alpha)$$
$$E_\alpha(-bt^\alpha)(D^\alpha y - by) = A_1 E_\alpha(at^\alpha)E_\alpha(-bt^\alpha)$$
$$D^\alpha(yE_\alpha(-bt^\alpha)) = \frac{A_1}{a-b}D^\alpha(E_\alpha(at^\alpha)E_\alpha(-bt^\alpha))$$

On integrating both side we get that is applying $D^{-\alpha}$ on both sides of above, we get

$$yE_\alpha(-bt^\alpha) = \frac{A_1}{a-b}(E_\alpha(at^\alpha)E_\alpha(-bt^\alpha)) + B$$

$$y = AE_\alpha(at^\alpha) + BE_\alpha(bt^\alpha) \quad \text{where} \quad A = \frac{A_1}{a-b}.$$

Therefore $y = AE_\alpha(at^\alpha) + BE_\alpha(bt^\alpha)$ is a solution of the differential equation. Thus we can state the following theorem

**Theorem 2:** The fractional differential equation

$$(D^\alpha - a)(D^\alpha - b)y(t) = 0$$

has solution of the form

$$y = AE_\alpha(at^\alpha) + BE_\alpha(bt^\alpha)$$

where A and B are constants.

Proof of the theorem is follows from the previous arguments.

Similarly one can generalized the solution of the differential equation in the following form

If

$$(D^\alpha - a_1)(D^\alpha - a_2)(D^\alpha - a_3)...(D^\alpha - a_n)y(t) = 0$$

with all $a_i$'s are distinct be a fractional differential equation with $0 \leq \alpha < 1$ then solution of the differential equation will be

$$y = \sum_{i=1}^{n} A_i E_\alpha(a_i t^\alpha)$$

where $A_i$ are arbitrary constants and $E_\alpha(a_i t^\alpha)$ is one parameter Mittag-Leffler function.

Let us consider the function

$$y = (At^\alpha + B)E_\alpha(at^\alpha)$$

Where A and B are constants.

Then

$$D^{\alpha} y = \Gamma(1+\alpha) A E_{\alpha}(at^{\alpha}) + (At^{\alpha} + B) a E_{\alpha}(at^{\alpha})$$
$$D^{\alpha} y - ay = \Gamma(1+\alpha) A E_{\alpha}(at^{\alpha}) + (At^{\alpha} + B) a E_{\alpha}(at^{\alpha}) - a\{(At^{\alpha} + B) E_{\alpha}(at^{\alpha})\}$$
$$= \Gamma(1+\alpha) A E_{\alpha}(at^{\alpha})$$
$$D^{2\alpha} y - a D^{\alpha} y = a \Gamma(1+\alpha) A E_{\alpha}(at^{\alpha})$$
$$= a(D^{\alpha} y - ay)$$
$$D^{2\alpha} y - 2a D^{\alpha} y + a^2 y = 0$$

Thus solution of the differential equation

$$D^{2\alpha} y - 2a D^{\alpha} y + a^2 y = 0$$

is

$$y = (At^{\alpha} + B) E_{\alpha}(at^{\alpha})$$

*A* and *B* are constants.

On the other hand consider the differential equation

$$(D^{\alpha} - a)^2 y = 0 \quad \text{or} \quad (D^{2\alpha} - 2a D^{\alpha} + a^2) y = 0$$
$$D^{2\alpha} y - 2a D^{\alpha} y + a^2 y = 0$$
(8)

Let $(D^{\alpha} - a) y = v$ then equation (8) reduce to the form

$$(D^{\alpha} - a) v = 0 \tag{9}$$

Solution of this differential equation is $v(t) = A_1 E_{\alpha}(at^{\alpha})$

$$(D^{\alpha} - a) y = A_1 E_{\alpha}(at^{\alpha})$$
$$E_{\alpha}(-at^{\alpha})(D^{\alpha} y - ay) = A_1 E_{\alpha}(at^{\alpha}) E_{\alpha}(-at^{\alpha})$$
$$D^{\alpha}[y E_{\alpha}(-at^{\alpha})] = A_1 = D^{\alpha}\left[\frac{A_1 t^{\alpha}}{\Gamma(1+\alpha)}\right]$$
$$y E_{\alpha}(-at^{\alpha}) = \frac{A_1}{\Gamma(1+\alpha)} + B$$
$$y = (At^{\alpha} + B) E_{\alpha}(at^{\alpha}) \quad \text{where} \quad A = \frac{A_1}{\Gamma(1+\alpha)}$$

*A* and *B* are constants.

Thus the following theorem can be stated

**Theorem 3:** The fractional differential equation

$$D^{2\alpha} y - 2aD^{\alpha} y + a^2 y = 0$$

has solution of the form

$$y = (At^{\alpha} + B)E_{\alpha}(at^{\alpha})$$

where $A$ and $B$ are constants.

The proof of the theorem is already explained in the previous arguments.

**Theorem 4:** Solution of the fractional differential equation

$$D^{2\alpha} y - 2aD^{\alpha} y + (a^2 + b^{\alpha}) y = 0$$

is of the form

$$y = E_{\alpha}(at^{\alpha})[A\cos_{\alpha}(bt^{\alpha}) + B\sin_{\alpha}(bt^{\alpha})].$$

Proof: The given differential equation can be written in the following form

$$((D^{\alpha} - a)^2 + b^2) y = 0 \quad \text{or} \quad (D^{\alpha} - a + ib)(D^{\alpha} - a - ib) y = 0 \tag{10}$$

Using theorem 3 we get the solution of the fractional differential (10) can be written in the following form

$$y = A_1 E_{\alpha}((a+ib)t^{\alpha}) + B_1 E_{\alpha}((a-ib)t^{\alpha})$$

$$y = A_1 E_{\alpha}(at^{\alpha}) E_{\alpha}(ibt^{\alpha}) + B_1 E_{\alpha}(at^{\alpha}) E_{\alpha}(-ibt^{\alpha})$$
$$= E_{\alpha}(at^{\alpha})[A_1\{\cos_{\alpha}(bt^{\alpha}) + i\sin_{\alpha}(bt^{\alpha})\} + B_1\{\cos_{\alpha}(bt^{\alpha}) - i\sin_{\alpha}(bt^{\alpha})\}]$$
$$= E_{\alpha}(at^{\alpha})[A\cos_{\alpha}(bt^{\alpha}) + iB\sin_{\alpha}(bt^{\alpha})]$$

Where $\quad A = A_1 + B_1 \quad$ and $\quad B = A_1 - B_1.$

Thus we get useful results.

## 6.0 Conclusions

There are several methods to solve fractional differential equations, and the solution depends on the type of fractional derivative used. Here we develop an analytical method to find the solutions of linear fractional differential equation, composed by Jumarie fractional derivative in terms of one parameter Mittag-Leffler function. Some well known properties of Mittag-Leffler have been used to find solution of the fractional differential equations. The solutions obtained are similar as the solutions obtained usual calculus obtained in terms the exponential function. This conjugation with ordinary calculus when Jumarie type fractional derivative is used to compose the fractional differential equations is useful in several physical problems.

# 7.0 References

[1] de Oliveira, E. C., J. A. T. Machado. A review of definitions of fractional derivatives and Integral. Mathematical Problems in Engineering. Hindawi Publishing Corporation. 2014. 1-6.

[2] Jumarie, G. Modified Riemann-Liouville derivative and fractional Taylor series of non-differentiable functions Further results, Computers and Mathematics with Applications, 2006. (51), 1367-1376.

[3] G. Jumarie, "On the solution of the stochastic differential equation of exponential growth driven by fractional Brownian motion," *AppliedMathematics Letters*, vol. 18, no. 7, pp. 817–826, 2005.

[4] G. Jumarie, "An approach to differential geometry of fractional order via modified Riemann-Liouville derivative," *Acta Mathematica Sinica*, vol. 28, no. 9, pp. 1741–1768, 2012.

[5] G. Jumarie, "On the derivative chain-rules in fractional calculus via fractional difference and their application to systems modelling," *Central European Journal of Physics*, vol. 11, no. 6, pp. 617–633, 2013.

[6] X. J. Yang, *Local Fractional Functional Analysis and Its Applications*, Asian Academic Publisher Limited, Hong Kong, 2011.

[7] X. J. Yang, *Advanced Local Fractional Calculus and Its Applications*, World Science, New York, NY, USA, 2012.

**[8]** Miller KS, Ross B. An Introduction to the Fractional Calculus and Fractional Differential Equations.John Wiley & Sons, New York, NY, USA; 1993.

[9 ] Das. S. Functional Fractional Calculus $2^{nd}$ Edition, Springer-Verlag 2011.

[10] Podlubny I. Fractional Differential Equations, Mathematics in Science and Engineering, Academic Press, San Diego, Calif, USA. 1999;198.

[ 11] M. Caputo, "Linear models of dissipation whose *q* is almost frequency independent-ii," *Geophysical Journal of the Royal Astronomical Society*, , 1967. vol. 13, no. 5, pp. 529–539.

[12] K M Kolwankar and A D. Gangal. Local fractional Fokker plank equation, Phys Rev Lett. 80 1998.

[13 ] Abhay Parvate, A. D. Gangal. Calculus of fractals subset of real line: formulation-1; World Scientific, Fractals Vol. 17, 2009.

[14] Abhay Parvate, Seema satin and A.D.Gangal. Calculus on a fractal curve in $R^n$ arXiv:00906 oo76v1 3.6.2009; also in Pramana-J-Phys.

[15] Abhay parvate, A. D. Gangal. Fractal differential equation and fractal time dynamic systems, Pramana-J-Phys, Vol 64, No. 3, 2005 pp 389-409.

[16]. E. Satin, Abhay Parvate, A. D. Gangal. Fokker-Plank Equation on Fractal Curves, Seema, Chaos Solitons & Fractals-52 2013, pp 30-35.

[17]. Erdelyi A. Asymptotic expansions, Dover (1954).

[18]. Erdelyi.A. (Ed). Tables of Integral Transforms. vol. 1, McGraw-Hill, 1954.

[19]. Erdelyi.A. On some functional transformation Univ Potitec Torino 1950.

**Acknowledgement**

Acknowledgments are to **Board of Research in Nuclear Science** (BRNS), Department of Atomic Energy Government of India for financial assistance received through BRNS research project no. 37(3)/14/46/2014-BRNS with BSC BRNS, title "Characterization of unreachable (Holderian) functions via Local Fractional Derivative and Deviation Function.